\documentclass[12pt]{amsart}
\usepackage{eucal}
\usepackage{amssymb}
\usepackage{mathrsfs}
\usepackage{upref}
\newif \ifwide
\newif \ifavnermargin
\def \makemargins{
\ifwide
        \headheight=6pt
        \oddsidemargin .25in
        \evensidemargin .25in
        \textwidth 6.00in
\else
\fi
\ifavnermargin
        \headheight=7pt
        \textheight=574pt
        \textwidth=432pt
        \topmargin=14pt
        \oddsidemargin=18pt
        \evensidemargin=18pt
\else 
\fi
}

\widetrue
\makemargins

\theoremstyle{plain}
\newtheorem{theorem}[subsection]{Theorem}
\newtheorem{proposition}[subsection]{Proposition}
\newtheorem{lemma}[subsection]{Lemma}
\newtheorem{corollary}[subsection]{Corollary}
\newtheorem*{mainresult}{Theorem \ref{main.result}}
\newtheorem*{mainresult.weighted}{Theorem \ref{main.weighted}}

\theoremstyle{definition}
\newtheorem{definition}[subsection]{Definition}

\theoremstyle{remark}
\newtheorem{remark}[subsection]{Remark}
\newtheorem*{ack}{Acknowledgements}

\newcount\timehh\newcount\timemm
\timehh=\time 
\divide\timehh by 60 \timemm=\time
\newcommand{\draftauthor}[1]{\author{#1
    {
      --- \protect \protect\sc\today\ ---
      \ifnum\timehh<10 0\fi\number\timehh\,:\,\ifnum\timemm<10 0\fi\number\timemm
      \protect \, \, \protect \bf DRAFT
    }
  }
}
\newcommand{\HH}{{\mathfrak H}}
\newcommand{\PP}{{\mathbb P}}

\newcommand{\CC}{{\mathbb C}}
\newcommand{\ZZ}{{\mathbb Z}}
\newcommand{\QQ}{{\mathbb Q}}

\newcommand{\ee}{{\mathrm{e}}}
\newcommand{\ii}{{\mathrm{i}}}

\newcommand{\TTT}{{\mathscr{T}}}

\newcommand{\MMM}{{\mathscr{M}}}

\renewcommand{\theta}{\vartheta}
\renewcommand{\tilde}{\widetilde}

\renewcommand{\mod}{\bmod}

\DeclareMathOperator{\Gr}{Gr}

\begin{document}

\title{Elliptic functions and equations of modular curves}

\newif \ifdraft
\def \makeauthor{
\ifdraft
        \draftauthor{Lev A. Borisov, Paul E. Gunnells and Sorin Popescu}
\else
\author{Lev A. Borisov}
\address{Department of Mathematics\\
Columbia University\\
New York, NY  10027}
\email{lborisov@math.columbia.edu}
\urladdr{http://www.math.columbia.edu/\~{}lborisov}
\author{Paul E. Gunnells}
\address{Department of Mathematics and Computer Science\\
Rutgers University\\
Newark, NJ  07102--1811
}
\email{gunnells@andromeda.rutgers.edu}
\author{Sorin Popescu}
\address{Department of Mathematics\\
Columbia University\\
New York, NY  10027\\
State University of New York at Stony Brook\\
Stony Brook, NY 11794-3651}
\email{psorin@math.columbia.edu, sorin@math.sunysb.edu}
\urladdr{http://www.math.columbia.edu/\~{}psorin/,
http://www.math.sunysb.edu/\~{}sorin/}
\fi
}

\thanks{The third author was partially supported by NSF grant DMS-9610205.}

\draftfalse
\makeauthor

\ifdraft
        \date{\today}
\else
        \date{May 22, 2000}
\fi


%
%

\begin{abstract}
Let $p\ge 5$ be a prime. We show that the space of weight one
Eisenstein series defines an embedding into $\PP^{(p-3)/2}$ of the
modular curve $X_1(p)$ for the congruence group $\Gamma_1(p)$ that is
scheme-theoretically cut out by explicit quadratic equations.
\end{abstract}
\maketitle

%
%
\section{Introduction}\label{introduction}
Modular curves are compactifications of quotients of the upper half
plane $\HH$ by arithmetic groups $\Gamma\subset {\mathrm
{SL}}_2(\QQ)$ acting on $\HH$ via M\"obius transformations. The group
varies depending on the type of polarization and, when appropriate,
the choice of the level structure (see \cite{shimura} for a full
account). Since $\Gamma$ is arithmetic, in other words commensurable
with ${\mathrm {SL}}_2(\ZZ)$, the quotient $\Gamma \backslash \HH$ is an
open Riemann surface that can be compactified to the modular curve
$X(\Gamma)$ by the addition of finitely many points (cusps).
 
Historically, several cases are of particular interest. The classical
one is where $\Gamma$ is the principal congruence subgroup
$\Gamma(p)\subset {\mathrm {SL}}_2(\ZZ)$ of level $p$, consisting of
all matrices congruent to the identity matrix modulo an odd prime
number $p$, and where the modular curve $X(p):=X(\Gamma(p))$
parameterizes elliptic curves with canonical level structure of order
$p$ (i.e. with fixed isomorphisms between their group of $p$-torsion
points with the Weil pairing, and the abstract group $(\ZZ/p\ZZ)^2$
with a certain standard symplectic paring). It is standard knowledge
(see for instance \cite{shimura}, Section 1.6) that $X(1)=\PP^1$,
having one cusp, and that $X(p)$ has genus $(p^2-1)(p-6)/24+1$ and
$(p^2-1)/2$ cusps if $p\geq 3$. Several natural embedded models of
$X(p)$ have been extensively studied starting perhaps with Felix Klein
(\cite{Klein1}, \cite{Klein2}). The best understood model is the
so-called $z$-modular curve in $\PP^{(p-3)/2}=\PP(V_-)$, which is
essentially the image in the negative eigenspace of the Heisenberg
involution, acting on $\PP^{p-1}$, of the $0$-section of the universal
elliptic curve over $X(p)$ under a morphism whose restriction to each
smooth fiber is induced by $p$ times the origin of that elliptic
curve. It is known that the $z$-modular curve is smooth \cite{velu},
Th\'eor\`eme 10.6.  The large symmetry group ${\mathrm
{PSL}}_2(\ZZ/p\ZZ)$ of the $z$-modular curve accounts for interesting
geometry and highly structured defining equations. When $p=7$, the
$z$-modular curve in $\PP^2$ is defined by Klein's quartic
$x_0^3x_1+x_1^3x_2+x_2^3x_0=0$, whose automorphisms group ${\mathrm
{PSL}}_2(\ZZ/7\ZZ)$ has maximal possible order for its genus
\cite{Klein1}, \cite{velu} (see also \cite{gross.popescu.1}, Example
2.10). For $p=11$, the $z$-modular curve in $\PP^4$ is defined by the
$4\times 4$-minors of the Hessian of the unique ${\mathrm
{PSL}}_2(\ZZ/11\ZZ)$-invariant cubic hypersurface $W\subset\PP^4$
given by the equation $v^2w+w^2x+x^2y+y^2z+z^2v = 0$. We refer for
these and other beautiful facts about the geometry of $X(11)$ to
\cite{Klein2}, pp. 153-156, \cite{edge}, \cite{adler.ramanan},
\cite{gross.popescu.2}. In general, one has the following beautiful
picture (cf. \cite{adler.ramanan}, Theorem 19.7, p. 56): There exists
a ${\mathrm {PSL}}_2(\ZZ/p\ZZ)$-equivariant isomorphism $\Phi:S^2(V_-)
\cong \Lambda^2(V_+)$, where $V_\pm$ denote the $\pm$-eigenspaces of
the Heisenberg involution, such that by means of $\Phi$ the $z$-curve
coincides with ${\nu_2}^{-1}(\Gr(2,V_+)),$ where
$\Gr(2,V_+)\subset\PP(\Lambda^2(V_+))$ denotes the Pl\"ucker embedding
of the Grassmannian of 2-dimensional linear subspaces in $V_+$, and
$\nu_2:\PP(V_-)\to \PP(S^2(V_-))$ denotes the quadratic Veronese
embedding. See also \cite{gross.popescu.1}, Corollary 2.9 and Example
2.10, and \cite{gross.popescu.2}, \S 2, and proof and discussion
before of Lemma 2.1, for a different approach to this description of
the $z$-curve. It is conjectured that the $z$-modular curve is
linearly normal and that the rank $2$ vector bundle corresponding to
the embedding into $\Gr(2,V_+)$ is stable (see \cite{adler.ramanan},
\cite{dolgachev} for details).

In this paper we investigate natural equations for
$X_1(p):=X(\Gamma_1(p))$, the modular curve for the congruence
subgroup $\Gamma_1(p)$ of matrices satisfying
$$\Bigl\{
\Bigl(\begin{array}{cc}
a&b\\
c&d
\end{array}\Bigr)\in{\mathrm {SL}}_2(\ZZ)\Bigm|
\Bigl(\begin{array}{cc}
a&b\\
c&d
\end{array}\Bigr)=
\Bigl(\begin{array}{cc}
1&*\\
0&1
\end{array}\Bigr)\mod p\Bigr\},
$$
where $p\ge 3$ is a prime. Our point of view is quite different from
the one described above for $X(p)$---our embedding uses certain modular forms
of weight one, whereas the embedding of the $z$-modular curve uses
forms of weight $(p-3)/(2p-12)$.  It would be very interesting to
somehow relate these two approaches.

In our case the modular curve $X_1(p)$ parameterizes elliptic curves
with a choice of non-trivial $p$-torsion
point. Furthermore, $X_1(p)$ is a smooth curve of genus
$(p-5)(p-7)/24$ with $(p-1)$ cusps for $p\geq 5$ (see \cite{shimura}). 
We consider in this paper the linear system on
$X_1(p)$ defined by weight one Eisenstein series, and show 
that it defines an embedding of $X_1(p)$ into $\PP^{(p-3)/2}$ (see
Corollary \ref{modular.embed}).  Our main result is an
explicit description of the equations defining this embedding:

\begin{mainresult}
Let $p\geq 5$ be a prime number.  The space of weight one
Eisenstein series defines an embedding of the modular curve $X_1(p)
\subset\PP^{(p-3)/2}$ that is scheme-theoretically cut out by 
the quadratic equations
\begin{align*}
(p-4)(s_as_b+s_bs_c+s_cs_a)=&
2s_a^2+2s_b^2+2s_c^2 
- \frac 4{p-2}\sum_{k\neq 0}s_k^2\\
&+\sum_{k\neq 0,a}s_ks_{a-k}+
\sum_{k\neq 0,b}s_ks_{b-k}+\sum_{k\neq 0,c}s_ks_{c-k},
\end{align*}
for all $a,b,c \in(\ZZ/p\ZZ)^*$ with $a+b+c=0$, where
${\{s_a\}}_{a\in(\ZZ/p\ZZ)^*}$ is a suitable system of coordinates 
with $s_{-a}=-s_a$.
\end{mainresult}

In contrast with the conjectural picture for $X(p)$, our embedding
$X_1(p)\subset\PP^{(p-3)/2}$ is in general neither linearly nor
quadratically normal. In fact, the failure of linear or
quadratic normality reflects the existence of certain modular
forms (see Remark \ref{fail} for the precise statement).

The equations of $X_1(p)$ can be greatly simplified by considering a
related embedding of $X_1(p)$ into the weighted projective space
$\PP(1,\ldots,1,2,\ldots,2)$ with $(p-1)/2$ variables $s_a$ of weight
one and symmetry as above, and with $(p-1)/2$ variables $t_a$ of
weight two and symmetry $t_{-a}=t_a$:

\begin{mainresult.weighted}
The quadratic relations 
$$
s_as_b+s_bs_c +s_cs_a + t_a + t_b + t_c = 0,
$$ 
for all $a,b,c \in(\ZZ/p\ZZ)^*$ with $a+b+c=0$, scheme-theoretically
cut out the modular curve $X_1(p)\subset \PP(1,\ldots,1,2,\ldots,2)$.
\end{mainresult.weighted}

The paper is structured as follows. In Section \ref{embedding} we
recall first results from \cite{bor.gunn.1} and \cite{bor.gunn.2}
concerning modular forms for $\Gamma_1(p)$, then show that weight one
Eisenstein series define an embedding of the modular curve and compute
equations vanishing on the image of this embedding (Proposition
\ref{equations}).  In Section \ref{diff} we introduce a system of
differential equations \eqref{diff.eq} that mimics the system
satisfied by elliptic functions with poles of order one along a
subgroup of order $p$.  We construct Laurent series solutions to this
system, and show that they satisfy the defining quadratic equations of
a (possibly degenerate) elliptic normal curve in $\PP^{p-1}$. Finally, in
Section \ref{modular.curve} we use a deformation argument to show that
a deformation of the solutions of the system of differential equations
leads to a deformation of elliptic curve, thus eventually showing that
the equations in Theorem \ref{main.result} cut out
scheme-theoretically the modular curve.

It is worth mentioning that the elliptic functions appearing in this paper
can be interpreted as discrete versions of the solutions of the
associative Yang-Baxter equation that was introduced in a recent paper of 
Polishchuk \cite{pol}.

\begin{ack} It is a pleasure to thank Igor Krichever and Pavel Etingof
for useful discussions. 
We are also grateful to Dan Grayson and Mike Stillman for
{\sl Macaulay2} \cite{macaulay}, which helped us to carefully check
the structure of the equations described in this paper.
\end{ack}

%

%

\section{Embedding modular curves by Eisenstein series of weight one}
\label{embedding}
Let $p$ be a natural number, and let $\Gamma_1(p)$ be as in Section \ref{introduction}.
For any positive integer $r$, let $\MMM_r(p):=\MMM_{r} (\Gamma _{1} (p))$ be
the
$\CC$-vector space of weight $r$ holomorphic modular forms for 
$\Gamma_{1}(p)$. 

Most results of this section hold with minor modifications for
non-prime levels, but for simplicity we consider only odd prime levels $p$.
Let $$\theta(z,\tau)=q^{\frac{1}{8}} (2 \sin \pi z)
\prod_{l=1}^{\infty}(1-q^l) \prod_{l=1}^{\infty}(1-q^l \ee^{2 \pi \ii
z})(1-q^l \ee^{-2 \pi \ii z}),$$ be the Jacobi theta function, where
$q=\ee^{2 \pi \ii \tau}$ (cf. \cite{Chandra}). As a theta function
with characteristic this is $\theta_{1,1}(z,\tau)$
(cf. \cite{Mumford}).  The following functions on the upper half plane
$$s_a(\tau):=\frac1{2\pi\ii}\frac d{dz} \ln \theta(z,\tau)|_{z=a/p}=
\frac{\xi^a+1}{2(\xi^{a}-1)} - \sum _{d}q^{d}\sum _{k|d} (\xi^{ka}-
\xi^{-ka}),
$$
where $\xi=\ee^{2\pi\ii/p}$ and $a$ is an integer with $a\neq 0 \mod p$,  
are modular forms of weight one with respect to the group 
$\Gamma_1(p)$ (see for instance \cite{Lang}). 
In fact, $s_a=s_{a+p}$ and thus we will regard 
the subscripts $a$ as elements of $(\ZZ/p\ZZ)^{*}$. Moreover
it follows directly from the definition that $s_{-a}=-s_{a}$.
Notice that each $s_a(\tau)$ is an Eisenstein series.
Our notation here
differs slightly from that of \cite{bor.gunn.1}, \cite{bor.gunn.2} 
and  \cite{bor.gunn.3}. What we denote by $s_{a}$ here was  
called $s_{a/p}$ in those papers. We will 
need the following results:

\begin{theorem}\label{getall}
\emph{(\cite{bor.gunn.1}, \cite{bor.gunn.2}, 
\cite{bor.gunn.3})}
The subring $\TTT_*(p)\subset \MMM_*(p)$ generated by the $s_a$ is 
Hecke-stable.  Moreover, for weights $k\geq 3$ we have $\TTT_k(p)=\MMM_k(p)$,
while for weight two $\TTT_2(p)$ is the direct sum of the space of
Eisenstein series and the span of Hecke eigenforms of analytic rank
zero.
\end{theorem}

\begin{corollary}\label{modular.embed}
The map $X_1(p)\to \PP^{(p-3)/2}$ that sends $\tau$ to
$\{s_a(\tau)\}$ defines a closed embedding of $X_1(p)$ into
$\PP^{(p-3)/2}$. 
\end{corollary}

\begin{remark}
One may prove directly this corollary without using the
calculations in \cite{bor.gunn.3}. For instance, to show that the map
defined by the weight one forms $s_a$ separates points, different from
the cusps, one reconstructs the values of $s_{a/p}^{(k)}$, in the
notation of \cite{bor.gunn.1}, up to a scaling factor $c^k$. This
allows reconstruction of the elliptic functions
$$\frac{\theta(a/p-cz,\tau)\theta_z(0,\tau)}
{\theta(-cz,\tau)\theta(a/p,\tau)}=
\exp(\sum_{k\geq 1}{(2\pi\ii z)}^k\frac{s_{a/p}^{(k)}}{k!}).
$$
Their poles and zeros determine uniquely up to scaling a lattice $\Lambda$ 
in $\CC$ together with an element of order $p$ in $\CC/\Lambda$.
Therefore, $\{s_a\}_{a\in(\ZZ/p\ZZ)^*}$ determine uniquely a 
point of $X_1(p)$. We leave the details to the reader.
\end{remark}

Since the $s_a$'s define an embedding, we will abuse
notation and denote in the sequel the image of the embedding of the
modular curve in $\PP^{(p-3)/2}$ simply by $X_1(p)$.\par\medskip

The following proposition follows easily from Lemma 4.8 in 
\cite{bor.gunn.1} for $N=2$, and Propositions 4.7 and 4.8
in  \cite{bor.gunn.2}:
\begin{proposition}\label{equations}
The weight two modular forms 
$$
t_a(\tau):=\frac{1}{2}(2\pi\ii)^{-2}(\frac{\theta_{zz}(a/p,\tau)}{\theta(a/p,\tau)}
-\frac{\theta_{zzz}(0,\tau)}{3\theta_{z}(0,\tau)})=
\frac 1{12} - \sum_{d>0}q^d\sum_{k|d}\frac dk
(\xi^{ak}+\xi^{-ak})
,
$$
where $a\in (\ZZ/p\ZZ)^*$, satisfy 
the following quadratic relations 
\begin{equation}\label{quad.rel}
s_as_b+s_bs_c +s_cs_a + t_a + t_b + t_c = 0
\end{equation}
for every triple of numbers  $a,b,c\in (\ZZ/p\ZZ)^*$ with $a+b+c=0$.
In particular, the symmetry of the $s_a$'s implies $t_{-a}=t_a$.
\end{proposition}


In fact the set of equations \eqref{quad.rel} allows one to express
$\{t_a\}$ in terms of $\{s_a\}$, and thus gives rise to degree two
relations involving only $\{s_a\}$. The goal of this paper is to show
that these relations cut out the curve $X_1(p)$ scheme-theoretically.

We will need the following results in the remaining two sections.
First of all, notice that the group $(\ZZ/p\ZZ)^*$ acts on $s_a$ and
$t_a$ by $s_a\rightsquigarrow s_{ak}$, $t_a\rightsquigarrow t_{ak}$.
Moreover, the action of the Fricke involution (see for example
\cite{Lang}) on $X_1(p)$ also lifts to an action $w_p$ on $\{s_a\}$
and $\{t_a\}$:

\begin{proposition}
The Fricke involution $w_p$ acts on the ring 
$\CC[s_1,\ldots, s_{p-1}, t_1,\ldots, t_{p-1}]$ as follows:
\begin{align*}
w_p(s_a)&={(-p)}^{-1/2}\sum_{k=1}^{p-1}\xi^{ka} s_k,\\
w_p(t_a)&={(-p)}^{-1}\sum_{k=1}^{p-1}\xi^{ka}(-s_k^2+2t_k).
\end{align*}
This transformation preserves the subscheme cut out by the quadratic
equations \eqref{quad.rel} and the symmetry relations on $s$ and $t$.
\end{proposition}

\begin{proof}
If $a+b+c=0$, where $a,b,c\neq 0$, then 
\begin{align*}
(-p)w_p (&s_as_b+s_bs_c+s_cs_a+t_a+t_b+t_c)=\\
&=\sum_{k,l\neq 0}(\xi^{ka+lb}+\xi^{-la+(k-l)b}+\xi^{(l-k)a-kb})s_ks_l+
\sum_{k\neq 0}(\xi^{ka}+\xi^{kb}+\xi^{kc})(-s_k^2+2t_k) \\
&=-\sum_{\substack{k,l\neq 0\\k\neq l}}\xi^{ka+lb}(s_{-k}s_l+s_{-k}s_{k-l}+
s_ls_{k-l}+t_{-k}+t_l+t_{k-l}).
\end{align*}
The fact that $w_p$ is indeed the Fricke involution follows from 
the details of the calculation in \cite{bor.gunn.1}.
\end{proof}

In what follows we shall also need a description of the cusps of 
$X_1(p)$.
\begin{proposition}
There are $p-1$ cusps on $X_1(p)$. They can all be obtained
by the action of $(\ZZ/p\ZZ)^*$ and $w_p$ from the point with coordinates
$$s_a = \lambda (1-\frac {2a}p),~1\leq a\leq p-1.$$
\end{proposition}

\begin{proof}
Every cusp on $X_1(p)$ can be obtained from $\ii\infty$ by 
means of $(\ZZ/p\ZZ)^*$ and the Fricke involution. It remains to
recall now the value of $s_a(\tau)$ at $\tau=0$, which is computed
by means of the Fricke involution.
\end{proof}

\section{Differential equations and quadratic relations}\label{diff}

We start by analyzing the system of quadratic equations
\eqref{quad.rel}, where $\{s_a,t_a\}$ will be now considered as free
variables rather than functions on the upper half-plane. To each
solution of this system we will associate a curve in $\PP^{p-1}$ that,
as we will see in Section \ref{modular.curve}, turns out to be a
$\ZZ/p\ZZ$-equivariant elliptic normal curve, or a degeneration
thereof. The main idea is to set up a suitable system of differential
equations that mimics the system satisfied by the elliptic functions
$$z\longmapsto \frac{\theta(a/p-z,\tau)\theta_z(0,\tau)}
{\theta(-z,\tau)\theta(a/p,\tau)}.
$$
We then construct Laurent series solutions to these equations, and show
that they satisfy certain quadratic relations (see \eqref{rel1}, \eqref{rel2}
bellow) that are the defining equations of an ``elliptic'' normal curve
in $\PP^{p-1}$.

\begin{definition}\label{ODE}
For each solution $\{s_a,t_a\}$ of the quadratic relations \eqref{quad.rel}
we introduce the system of ordinary differential equations 
for functions $r_a(z)$, 
$a\in (\ZZ/p\ZZ)^*$:
\begin{equation}\label{diff.eq}
\frac {d r_a}{dz} = -\frac {1}{p-2}
\bigl( \sum_{k\neq 0,a} r_kr_{a-k} + 2r_as_a\bigr).
\end{equation}
\end{definition}

We will be looking for solutions with a pole of order one at $z=0$:
\begin{definition}\label{standard} 
A solution $\{\hat r_a\}$, $a\in (\ZZ/p\ZZ)^*$ 
of the system \eqref{diff.eq} with Laurent series  expansion
$$\hat r_a(z) = \frac 1z + s_a + t_a z + \cdots$$
will be called a {\em standard solution} of this system of differential equations.
Such a solution will be interpreted as an element 
$\hat r=\sum_{a=0}^{p-1}\hat r_a\xi^a$
in $\CC[\xi]/(\xi^p-1)[[z,z^{-1}]]$, where now $\hat r_0$ is defined to be $1$.
\end{definition}

\begin{proposition}
Standard solutions exist and are unique. They are convergent in a small 
punctured disc around $z=0$ in the sense of convergence of functions 
with values in $\CC[\xi]/(\xi^p-1)$.
\end{proposition}

\begin{proof}
Observe first that the relations among $\{s_a,t_a\}$ imply that $\hat
r$ satisfies the differential equations up to first order in $z$.
More explicitly, by summing \eqref{quad.rel} up for all triples
$(k,a-k,-a)$ and using the fact that $\sum_{k\neq 0}s_k=0$, we obtain
the necessary relation
\begin{equation}
\label{elim.t}
-(p-2)t_a=\sum_{k\neq 0,a}s_ks_{a-k}+\sum_{k\neq 0,a}(t_k+t_a)+2s_a^2.
\end{equation}
On the other hand, the differential equations lead to recursive relations on the coefficients
$\hat r_{a,n}$ of $\hat r_a =\sum_n \hat r_{a,n}z^k$ of the form
$$
n\hat r_{a,n} = -
\frac {1}{p-2}
\bigl( \sum_{k\neq 0,a} \hat r_{k,n}+\hat r_{a-k,n} + 
\sum_{k\neq 0,a}\sum_{d=0}^{n-1}
\hat r_{k,d}\hat r_{a-k,n-1-d} + 2s_a \hat r_{a,n-1}
\bigr)
$$
Observe now that these linear equations define $\hat r_{a,n}$ uniquely. 
Moreover, it is straightforward to show by induction that $|\hat r_{a,n}| \leq c^n$ for some constant $c$. This follows since the inverse of the 
matrix defining the above recursion is $\frac 1n{\bf id} + O(1)$ as 
$n\to\infty$. Thus the convergence of the series follows, and this concludes the proof. 
We remark that the constant $c$ depends on $\{s_a,t_a\}$. 
\end{proof}

We will now introduce free variables $r_0, \ldots, r_{p-1}$ to study
the algebraic relations satisfied by the $\hat r_a(z)$. These free variables
should not be confused with the complex-valued functions in Definition 
\ref{ODE}.

\begin{definition}
We introduce now the following two sets of relations in the polynomial ring
$\CC[r_0,\ldots,r_{p-1}]$. 
\begin{equation}\label{rel1}
R_{a,b,c,d}:= r_a r_b - r_c r_d - r_0r_{a+b}(s_a+s_b-s_c-s_d),
\ a+b=c+d\neq 0,\ a,b,c,d\neq 0,
\end{equation}
and
\begin{equation}\label{rel2}
R_{a,b,-a,-b}:= r_a r_{-a} -r_b r_{-b} - r_0^2( -s_a^2 +2 t_a +s_b^2 -2t_b),~
a,b\neq 0.
\end{equation}
\end{definition}

These relations define an homogeneous ideal $I_{s,t}$ in the ring 
$\CC[r_0,r_1,\ldots,r_{p-1}]$. 
The goal of the rest of the section is to calculate the Hilbert function
of $\CC[r_0,r_1,\ldots,r_{p-1}]/I_{s,t}$.

\begin{proposition}\label{upper.bound}
For all $n>0$, 
$$\dim_\CC \bigl(\CC[r_0,r_1,\ldots,r_{p-1}]/I_{s,t} \bigr)_n \leq np.$$
The dimension of the zero graded component is one.
\end{proposition}

\begin{proof}
Observe that this ring is $\ZZ/p\ZZ$-graded by the subscripts of
$r_a$.  The quadratic relations \eqref{rel1} and \eqref{rel2} imply
that the dimension of each $\ZZ/p\ZZ$-graded component of
$\bigl(\CC[r_0,r_1,\ldots,r_{p-1}]/I_{s,t} \bigr)_n$ is at most
$n$. Indeed, a spanning set is constructed inductively by multiplying
the previous spanning set by $r_0$ and by adding any monomial of the
correct $\ZZ/p\ZZ$-weight that does not contain $r_0$. We need to show
that modulo $r_0$ any two monomials with the same $\ZZ/p\ZZ$-weight
are equivalent. We use induction on the degree and the number of equal
factors.  If one monomial is $r_{a_1}r_{a_2}\cdots$ and the other one
is $r_{b_1}r_{b_2}\cdots$, then we may reduce the first monomial to
$r_{b_1}r_{a_1+a_2-b_1}\cdots$ unless $b_1=a_1+a_2$. This argument
shows that the only difficulty could occur in the case of monomials
$r_a^n$ and $r_b^n$ with $2a=b$ and $2b=a$. However, this is not
possible for $p>3$.
\end{proof}

We will now show that $\hat r$ satisfies the quadratic relations 
\eqref{rel1} and \eqref{rel2}.
\begin{proposition}\label{key}
Define $\hat R_{a,b,c,d}(z)$, for $a+b=c+d,\ a,b,c,d\neq 0$ by 
$$
\hat R_{a,b,c,d}(z) = 
\hat r_a \hat r_b - \hat r_c \hat r_d - \hat r_0\hat r_{a+b}
(s_a+s_b-s_c-s_d),~a+b\neq 0,
$$
$$
\hat R_{a,b,-a,-b}(z) = 
\hat r_a \hat r_{-a} -\hat r_b \hat r_{-b} - 
\hat r_0^2( -s_a^2 +2 t_a +s_b^2 -2t_b).
$$
Then $\hat R_{a,b,c,d}(z) \equiv 0$.
\end{proposition}

\begin{proof}
The differential equations for the $\hat r$'s imply the following differential
equations for the $\hat R$'s:
\begin{align*}
\frac {d \hat R_{a,b,c,d}}{dz} =& -\frac 1{p-2} \bigl(
\sum_{k\neq 0,a,c,a+b} \hat r_k \hat R_{a-k,b,c-k,d} 
+
\sum_{k\neq 0,b,d,a+b} \hat r_k \hat R_{b-k,a,d-k,c}\\ 
&+\hat r_{a+b}(\hat R_{a,-a,c,-c} +\hat R_{b,-b,d,-d})
+(s_a+s_b+s_c+s_d)\hat R_{a,b,c,d}\bigr ),\\
\frac {d \hat R_{a,b,-a,-b}}{dz} =& -\frac 1{p-2}
\sum_{k\neq 0,a,b}\hat r_k (\hat R_{-a,a-k,-b,b-k}+\hat R_{a,-a+k,b,-b+k}).
\end{align*}
To derive the first set of differential equations, we use the formula
for the derivatives of $\hat r$ and then collect terms with the 
common factor $\hat r_k$ (where $k$ is the index of summation in
equation \eqref{diff.eq}). We then split off products of $\hat r_k$
and $\hat R$ to get
\begin{align*}
&(2-p)\frac {d \hat R_{a,b,c,d}}{dz} = 
\sum_{k\neq 0,a,c,a+b} \hat r_k \hat R_{a-k,b,c-k,d} +
\sum_{k\neq 0,b,d,a+b} \hat r_k \hat R_{b-k,a,d-k,c} 
\\
&
+
\sum_{k\neq 0,a,c,a+b} \hat r_k\hat r_{a+b-k}(s_{a-k}+s_{b} -s_{c-k} -s_d)
+
\sum_{k\neq 0,b,d,a+b} \hat r_k\hat r_{a+b-k}(s_{b-k}+s_{a} -s_{d-k} -s_c)
\\
&
-\sum_{k\neq 0,a+b}\hat r_k\hat r_{a+b-k}(s_a+s_b-s_c-s_d)
+\hat r_{a+b}(\hat r_a\hat r_{-a}+\hat r_b\hat r_{-b}
-\hat r_c\hat r_{-c}-\hat r_d\hat r_{-d})
\\&
+2\hat r_a\hat r_b(s_a+s_b)-2\hat r_c\hat r_d(s_c+s_d)
-2\hat r_{a+b}s_{a+b}(s_a+s_b-s_c-s_d).
\end{align*}
We observe that the sums that contain $\hat r_k\hat r_{a+b-k}$ would
cancel if the summation were over all $k\neq 0,a+b$, since
$k\rightsquigarrow a+b-k$ results in $a-k\rightsquigarrow k-b$,
$c-k\rightsquigarrow k-d$. Thus one may reduce all the terms to $\hat
r_{a+b}$ modulo $\hat R$.  The fact that $\{s_a,t_a\}$ satisfy
equations \eqref{quad.rel} insures that the coefficient of $\hat
r_{a+b}$ vanishes.  The differential equations for $\hat
R_{a,-a,b,-b}$ are treated similarly.

One can check that the equations on $\{s_a\}$ and $\{t_a\}$ imply that
$\hat R$ is zero to the order $z^0$. The vanishing of other Laurent
coefficients of $\hat R(z)$ is proved by induction on the degree of
$z$. It is easy to see that the recursion matrix is invertible
starting from the coefficients at $z^2$, because the diagonal entries
dominate the rows.  The case of the coefficients by $z^1$ is handled
separately as follows.

It is easy to see that $\hat r_a(z)=1/z + s_a +t_a z +u_a z^2 +\cdots $,
where 
$$
u_a = \frac 1{(p-3)} (\sum_{k\neq 0,a} s_{k-a}t_k -s_at_a).
$$
To show that the coefficient of  $\hat R_{a,b,c,d}$ at $z^1$ is zero,
it is enough to show
\begin{equation}
\label{u.rel}
s_at_b+s_bt_a -(s_a+s_b)t_{a+b} + u_a+u_b + 2u_{a+b}=0
\end{equation}
for all $a,b\neq 0$ such that $a+b\neq 0$. 
To prove this, for every $k\neq 0,-a,b$, we use the relations \eqref{quad.rel}
for triples $(-k-a,k,a)$, $(k-b,-k,b)$, and $(-k+b,k+a,-a-b)$ multiplied
by $(s_{k-b}+s_b)$, $(s_{k+a}+s_a)$, and $(s_a+s_b)$ respectively
to get 
\begin{align*}
&s_as_{a+b}(s_{k+a}-s_{k-b}) + s_as_b(s_{k-b}-s_{k+a})
+ s_bs_{a+b}(s_{k+a}-s_{k-b})
\\&
-(s_a+s_b)t_{a+b}+s_at_b+s_bt_a+
(s_a+s_b)t_k-s_at_{k+a}-s_bt_{k-b}
\\&+
s_{k-b}t_k+s_{k-b}t_{k+a}+s_{k-b}t_a-s_{k+a}t_k -s_{k+a}t_{k-b}-s_{k+a}t_b
=0.
\end{align*}
We then sum the above equations for all $k\neq 0,-a,b$. Finally we use $\sum_{k\neq 0}
s_k = 0$ and the relation \eqref{quad.rel} for $(a,b,-a-b)$ to 
get 
\begin{align*}
s_at_a+s_bt_b+2s_{a+b}t_{a+b}=&
(p-3)(s_at_b+s_bt_a-(s_a+s_b)t_{a+b})\\
&+\sum_{k\neq 0,b}s_{k-b}t_k+\sum_{k\neq 0,a}s_{k-a}t_k+
2\sum_{k\neq 0,a+b}s_{k-a-b}t_k,
\end{align*}
which is equation \eqref{u.rel} above.
To finish the proof, observe now that the differential equations for
$\hat R_{a,-a,b,-b}$
imply that it vanishes to the order $z^1$ as long as all 
$\hat R_{a,b,c,d}$ with
$a+b=c+d\neq 0$ vanish. Alternatively, observe that 
$\hat R_{a,-a,b,-b}(z)$ is even, because $\hat r_{a}(-z)=-\hat r_{-a}(z)$.
\end{proof}

\begin{theorem}
For all $n>0$, 
$$\dim_\CC \bigl(\CC[r_0,r_1,\ldots,r_{p-1}]/I_{s,t} \bigr)_n = np.$$
In other words, the Hilbert function of this ring is the same as the 
Hilbert function of the homogeneous coordinate of an elliptic normal
curve of degree $p$ in $\PP^{p-1}$.
\end{theorem}

\begin{proof}
By Proposition \ref{key} we can use $\hat r$ to map 
$\CC[r_0,r_1,\ldots,r_{p-1}]/I_{s,t}$ to $\CC[\xi]/(\xi^p-1)[[z,z^{-1}]]$. 
We observe that each $\ZZ/p\ZZ$-graded component of degree $n$
is mapped onto a space of dimension at least $n$, because polynomials
of degree $n$ in $\hat r$ and fixed $\ZZ/p\ZZ$-grading can have an 
arbitrary singular part in their Laurent expansions. The only exception
is the zero graded component, where the coefficient by $z^{-1}$ cannot
be chosen freely (in fact, it is always zero), but $\hat r_0^n$ allows
us to freely choose the constant term.
This provides a lower bound on the Hilbert
function, which together with Proposition \ref{upper.bound} finishes the proof.
\end{proof}

\section{The modular curve}\label{modular.curve}

As a first step to proving our main result, we show the following:

\begin{theorem}\label{set-theoretic}
The quadratic relations \eqref{quad.rel} on $s_a$ and $t_a$ cut out 
set-theoretically the modular
curve $X_1(p)\subset \PP(1,\ldots,1,2,\ldots,2)$.
\end{theorem}

\begin{proof}
Let $\{s_a,t_a\}$ be a solution to the system \eqref{quad.rel}. Then
the results of Section \ref{diff} allow us to define a dimension one subscheme
$C$ of $\PP^{p-1}$ cut out by the ideal $I_{s,t}$. 
Notice that points $p_k,~k=0,\dots ,p-1$ with coordinates 
$(w_k^1:w_k^2:\cdots:w_k^{p-1}:0)$, where $w_k=\exp(2\pi\ii k/p)$, lie on $C$.
Moreover, $\hat r$ defines a map from the disjoint union of the neighborhoods
of $p$ non-singular points in $\CC$ to $C$. As a result, the embedding
$$\CC[r_0,r_1,\ldots,r_{p-1}]/I_{s,t}\to \CC[\xi]/(\xi^p-1)[[z,z^{-1}]]$$
factors through  
$$\CC[r_0,r_1,\ldots,r_{p-1}]/I_{s,t}\to \CC[r_0,r_1,\ldots,r_{p-1}]/J_{s,t},$$
where $J_{s,t}$ is the defining ideal of the reduced irreducible components
of $C$ that pass through points $p_k$. 
Indeed, a map to $\CC[\xi]/(\xi^p-1)[[z,z^{-1}]]$ is determined by a 
collection of $p$ maps to the ring of usual Laurent series, which is an 
integral domain.
This shows that $I_{s,t}=J_{s,t}$. In particular, $C$ is reduced and has no 
isolated closed points. 

Suppose now that $C$ is singular. Let  $p=(\rho_0:\cdots:\rho_{p-1})$ be a singular point
on $C$. We will first consider the case $\rho_0\neq 0$, which we can
normalize to get $\rho_0=1$. Let 
$p_\epsilon=(1:\rho_1+\epsilon u_1:\cdots:\rho_{p-1}+\epsilon u_{p-1})$ be a 
$\CC[\epsilon]/\epsilon^2$-point of $C$
that reduces to $p$. Set up the differential equations \eqref{diff.eq} 
for power series $\tilde r_a(z)=(\rho_a+\epsilon u_a)+ O(z)$.
The argument of Proposition \ref{key} can be extended to the 
power series over the ring of dual numbers. Unless the solutions
$\tilde r_a$ are constant, this allows one to move a
$\CC[\epsilon]/\epsilon^2$-point to a nearby location on $C$, which 
is impossible since $C$ is reduced. Therefore, the solution $\tilde r_a$ must
be constant, which implies
\begin{equation}\label{noflow}
0=-\frac {1}{p-2} \bigl( \sum_{k\neq 0,a} \rho_k\rho_{a-k} + 2\rho_as_a\bigr).
\end{equation}

\begin{lemma}
If the point $p$ satisfies the equation \eqref{noflow}, then $\{s_a,t_a\}$
correspond to a cusp of $X_1(p)$.
\end{lemma}

\begin{proof}
For each $a,b,a+b\neq 0$ let 
$$x_{a+b}:=\rho_a\rho_b-\rho_{a+b}(s_a+s_b).$$
Observe that because of relation \eqref{rel1}, $x_{a+b}$ 
depends on the sum $a+b$ only.
Equation \eqref{noflow} then implies $x_{a}=0$ for all $a\neq 0$.
As a consequence, we have 
\begin{equation}\label{noflow2}
\rho_a\rho_b=\rho_{a+b}(s_a+s_b)
\end{equation}
for all $a,b,a+b\neq 0$. There are now two cases to consider.

{\em Case 1.} One of the $\rho_a$, $a\neq 0$ is zero. If any other $\rho_b$
is non-zero, then the above equation for a triple $b,b,2b\neq 0$
implies $\rho_{2b}\neq 0$. Then the above equation for $b,2b,3b$ gives
$\rho_{3b}\neq 0$ and so on. Eventually we get $\rho_a\neq 0$, which is a
contradiction.  Therefore $\rho_a=0$ for all $a\neq 0$.

Relations \eqref{rel2} now imply that for all $a,b\neq 0$,
$$
s_a^2-2t_a=s_b^2-2t_b.
$$
Let us denote by $x=s_a^2-2t_b$. Then relations \eqref{quad.rel} 
become 
$$(s_a+s_b+s_c)^2=3x$$
for all $a,b,c\neq 0, a+b+c=0$. We may scale $\{s_a\}$ so that
$3x=1$ to get 
$$
s_a+s_b+s_c = \pm 1,\qquad a,b,c\neq 0,\quad a+b+c=0.
$$
It is easy to see that for any given choice of signs, either there is no solution, or 
the solution is unique and rational.

The action of $(\ZZ/p\ZZ)^*$ allows us to assume that 
$s_1$ is the biggest of all $s_a$. We have
$2s_1-s_2=\pm 1$, therefore $s_1\leq 1$. Since $s_{p-1}=-s_1$,
$s_1$ must be positive, so $2s_1-s_2=1$ and $s_1\leq 1$.
Therefore $|s_a|\leq 1$ for all $a$. For every $a=2,\ldots,p-2$
we have 
$$s_1+s_a-s_{a+1}=\pm 1.$$
If the right hand side is $-1$, then we have $1+s_1=s_{a+1}-s_a$,
which can happen only if $s_1=1$, $s_{a+1}=1$, and $s_a=-1$.
This implies $2+s_{a-1}=\pm 1$, so $s_{a-1}=-1$. Then analogously
$s_{a-2}=-1$, and so eventually we get $s_1=-1$, a contradiction. Therefore,
we have $s_{a+1}=s_a+s_1-1$ for all $a=1,\ldots,p-2$. Together
with $s_{p-1}=-s_1$ this forces 
$$s_k = 1 - \frac {2k} p,$$
which means that $\{s_a\}$ has values corresponding to one of the cusps of 
$X_1(p)$.

{\em Case 2.} All $\rho_a$ are non-zero. Then all $s_a+s_b$ are 
non-zero for $a+b\neq 0$. Denote by $x_a=-\rho_a/\rho_{-a}$. Notice
that $x_ax_b=x_{a+b}$ for all $a,b$, when we set $x_0=1$. Therefore,
$x_a=w^a$ where $w$ is a $p$-th root of unity. We can multiply 
each $\rho_a$ by $w_1^a$ where $w_1$ is an appropriately chosen 
$p$-th root of unity to reduce ourselves to the case $x_a=1$.
So now we have $\rho_{-a}=-\rho_a$. For every $a+b+c=0$, $a,b,c\neq 0$
consider now the equations \eqref{noflow2} 
$$\rho_a\rho_b=-\rho_c(s_a+s_b),\quad \rho_c\rho_a=-\rho_b(s_c+s_a).$$
They imply $\rho_a^2=(s_a+s_b)(s_a+s_c)$, so $s_a s_b + s_b s_c + s_a s_c$
depends on $a$ only. Now the quadratic relations \eqref{quad.rel} imply that 
$t_a+t_b+t_c$ depends on $a$ only, if $a+b+c = 0\mod p$. This easily 
yields that all $t_k$ are the same. If we act on this set of
$\{s_k,t_k\}$ by the Fricke involution, see Section \ref{embedding}, we
get $\{s_k,t_k\}$ for which all $s_a^2-2t_a$ are the same. Then
the argument of case 1 finishes the proof of the lemma.
\end{proof}

\noindent
{\it Proof of Theorem \ref{set-theoretic} continued.}
Returning to the case of the singular point $p$ with $\rho_0=0$, observe
that the equations \eqref{rel1} and \eqref{rel2} imply that the point is
one of the points $p_k$. Similar arguments show then that it has to be
nonsingular.

Because of the equality $I_{s,t}=J_{s,t}$ and the presence of the $\ZZ/p\ZZ$ action, either $C$ is 
a smooth curve of degree $p$ and genus one, or it consists of $p$ 
lines. In the latter  case, the Hilbert function together with the symmetry 
forces the lines to form a $p$-gon. In particular it has singularities,
which implies that $\{s_a,t_a\}$ corresponds then to a cusp of the modular 
curve. 

So $C$ is an elliptic curve with a $\ZZ/p\ZZ$-action which permutes the points
$p_k$. Therefore, these points form a subgroup $S$ of order $p$ on $C$. 
Observe now that the embedding of $C$ by the $r_a$'s is given by the 
complete linear system $H^0({\mathcal O}_C(S))$.
Because the $r_a$'s are eigenvectors of the $\ZZ/p\ZZ$-action, 
we get 
$$\frac {r_a}{r_0} = \lambda_a
\frac{\theta(a/p-z,\tau)\theta_z(0,\tau)}
{\theta(-z,\tau)\theta(a/p,\tau)},
$$
where $\theta$ denotes the Jacobi theta function and $z$ is a uniformizing
parameter on the universal cover of the elliptic curve $C$ such that the 
points $p_k$ lift to $\frac kp+\ZZ$.  Since all the coefficients
of the products $r_ar_b$ in the relations \eqref{rel1} and \eqref{rel2}
are one, all $\lambda_a$ are equal to one.
Then the $s_a$'s can be uniquely 
determined from the equations \eqref{rel1} 
up to a multiplicative constant. More precisely, we rescale the $s_a$'s
and $r_a$'s to get $r_0=1$. Then for any two solutions $\{s_{a}'\}$ and
$\{s_{a}''\}$ their componentwise difference $\{\tilde s_a\}$ satisfies 
$\tilde s_a+\tilde s_b=\tilde s_c+\tilde s_d$ whenever $a+b=c+d\mod p$. This together with the 
symmetries of the $\tilde s_a$'s yields $\tilde s_a=0$. 
It remains to observe now that 
$s_a=\partial_z\log\theta(a/p,\tau)$ is a solution to \eqref{rel1}.
This finishes the proof.
\end{proof}

\begin{corollary}
The system of differential equations \eqref{diff.eq} has a solution in 
elliptic functions with poles of order one along a subgroup of order $p$
if and only if $\{s_a,t_a\}$ satisfy the quadratic relations 
\eqref{quad.rel} and do not correspond to a cusp.
\end{corollary}

Finally we prove now the main result of the paper:
\begin{theorem}\label{main.weighted}
The quadratic relations \eqref{quad.rel} on $s_a$ and $t_a$ cut out scheme-theoretically
the modular curve $X_1(p)\subset \PP(1,\ldots,1,2,\ldots,2)$.
\end{theorem}

\begin{proof}
We need to show that the scheme $\tilde X_1(p)\subset
\PP(1,\ldots,1,2,\ldots,2)$ cut out by these quadratic relations is
smooth. Let $\{s_a,t_a\}$ be a closed point of $\tilde X_1(p)$ that
is not a cusp. Let $\{s_a:t_a\}_\epsilon$ be a
$\CC[\epsilon]/\epsilon^2$-point of $\tilde X_1(p)$ that reduces to
$\{s_a:t_a\}$. Even though $s_a$ and $t_a$ are defined only up to
homothety (multiplication by $\lambda$ and $\lambda^2$, respectively), we
shall make a choice and fix a solution as a set of numbers.

Then the quadratic relations \eqref{rel1} and \eqref{rel2} define an
elliptic curve $C$ embedded into $\PP^{p-1}$ by a complete linear
system of degree $p$. A $\CC[\epsilon]/\epsilon^2$-point defines a
deformation of the homogeneous coordinate ring of this embedding. In
fact, one can set up the system of differential equations on $\hat
r_a(z)$ that will now take values in $\CC[\epsilon]/\epsilon^2$
instead of $\CC$. Then its solutions will satisfy the quadratic
relations and will define a map from
$(\CC[\epsilon]/\epsilon^2)[r_0,\ldots,r_{p-1}]/I_{s,t}$ to the ring of
Laurent series with coefficients in the dual numbers. As before one shows
that the dimension of
$((\CC[\epsilon]/\epsilon^2)[r_0,\ldots,r_{p-1}]/I_{s,t})_n$ over $\CC$ is
$2pn$ and that it is a free module over $\CC[\epsilon]/\epsilon^2$
with a monomial basis as in Proposition \ref{upper.bound}.

Therefore we have a $\ZZ/p\ZZ$-equivariant deformation of
the embedding of an elliptic normal curve in ${\PP}^{p-1}$.
Moreover the $p$-torsion points 
$x_k$, $k=0,\ldots,p-1$ with coordinates 
$(r_0:\cdots:r_{p-1})=(0:\xi^k:\xi^{2k}:\cdots:\xi^{(p-1)k})$, 
where $\xi=\exp(2\pi\ii/p)$ are fixed under this deformation.

Such $\ZZ/p\ZZ$-equivariant embedded deformations are parameterized by
elements of
$$
{H^0(N_{C\mid{\PP}^{p-1}}
\otimes{\mathcal I}_{\{x_0,\ldots, x_{p-1}\}})}^{\ZZ/p\ZZ}=
H^0(N_{C\mid{\PP}^{p-1}}(-1))^{\ZZ/p\ZZ},$$
where  $N_{C\mid{\PP}^{p-1}}$ is the normal bundle of $C\subset\PP^{p-1}$.
An easy calculation shows that this space is 2-dimensional. Therefore
all such deformations come either from the scaling of the $s_a$'s and $t_a$'s,
or from deformations along the modular curve. This shows that the
scheme cut out by the quadratic equations in \eqref{quad.rel} is
nonsingular except possibly at the cusps. 

Thus it remains to show that in the neighborhood of the cusps the quadratic
relations \eqref{quad.rel} cut out a smooth curve. The action of $(\ZZ/p\ZZ)^*$
and the Fricke involution allows one to consider only the case 
$$s_a=1-\frac {2a}p,\quad t_a=s_a^2-\frac 13.$$
It suffices to calculate  the dimension of the tangent space at this point.
Denote the coordinates in the tangent space at this point in $\CC^{p-1}$
by $ds_a,dt_a,~a=1,\ldots,p-1$ with $ds_{p-a}=-ds_a$ and $dt_{-a}=d_{t_a}$.
Equations \eqref{quad.rel} yield in this tangent space
$$
(ds_a+ds_b+ds_c)=(d\hat t_a+d\hat t_b+d\hat t_c)
$$
for $a,b,c\in\{1,\ldots,p-1\}$ with $a+b+c=p$, where
$d\hat t_k=-dt_k+s_k ds_k$. 
 Two solutions are easy to 
see, namely 
$$ds_k=s_k,\quad dt_k=2t_k$$
and 
$$ds_k = \delta_{k}^1-\delta_{k}^{p-1},\quad d\hat t_k=2\delta_{k}^1+
2\delta_k^{p-1},$$
where $\delta$ is the Kronecker delta function.
By subtracting a linear combination of these solutions one can reduce any 
other solution to the case $ds_1=d\hat t_1=0$.
For every $a=1,\ldots,p-2$ we get
$$
ds_a-ds_{a+1} = d\hat t_a+d\hat t_{a+1}
$$
which implies $ds_2=-d\hat t_2$ and 
$$ds_k = -2d\hat t_2-2d\hat t_3 -\cdots- 2d\hat t_{k-1}-d\hat t_k,
\quad k=3,\ldots,p-1.$$
On the other hand, for every $a=1,\ldots,p-3$ 
one has 
$$ds_2+ds_a-ds_{a+2}=d\hat t_2+d\hat t_a+d\hat t_{a+2}$$
which now implies 
$$d\hat t_2=d\hat t_{a+1}.$$
Because $ds_{p-1}=0$, 
we deduce that  $d\hat t_2=0$, and then that all $ds_k$ and 
$dt_k$ must be zero. Thus the dimension of the tangent space is one, 
after we mod out by the rescaling factor.
\end{proof}

\begin{corollary}\label{main.result}
The embedding of $X_1(p)$ into $\PP^{(p-3)/2}$ induced by the weight
one modular forms $\{s_a\}$ is scheme-theoretically cut out 
by quadrics obtained by eliminating the $t$'s from the equations 
\emph{(\ref{quad.rel})}. Explicitly, $X_1(p)$ is cut out by
the quadratic relations
\begin{align*}
(p-4)(s_as_b+s_bs_c+s_cs_a)=&
2s_a^2+2s_b^2+2s_c^2 
- \frac 4{p-2}\sum_{k\neq 0}s_k^2\\
&+\sum_{k\neq 0,a}s_ks_{a-k}+
\sum_{k\neq 0,b}s_ks_{b-k}+\sum_{k\neq 0,c}s_ks_{c-k},
\end{align*}
for all $a,b,c \in(\ZZ/p\ZZ)^*$ with $a+b+c=0$.
\end{corollary}

\begin{proof}
This is immediate by using relations \eqref{elim.t} to eliminate
all $t_a$'s in \eqref{quad.rel}. Details are left to the reader.
\end{proof}

\begin{remark}\label{fail}
1) The embedding $X_1(p)$ into $\PP^{(p-3)/2}$ is generally not
linearly normal. This is due to the existence of weight one cusp
forms, and such forms are never in the span of $\{s_a(\tau)\}$. This first
happens at $p=23$, see \cite{stark}. Moreover, $X_1(23)$ is a
projection of a canonical curve.
\hfill\par 2) Furthermore, the embedding is not quadratically normal
in general, due to the existence of Hecke eigenforms of positive
analytic rank, see \cite{bor.gunn.2}. The first occurrence is at
$p=37$. The stability of the bundle of rank $(p-3)/2$, 
which is the kernel of the natural evaluation map, would imply 
through the known Koszul techniques a bound on the number of Hecke
eigenforms of non-zero analytic rank, which is linear in $p$.
\hfill\par 3) The quadratic relations \eqref{quad.rel} do not always
generate the ideal of the embedding in degrees three or more. Macaulay
\cite{macaulay} calculations show that this first happens at $p=43$.
\end{remark}

%
%
%
%

%
%

%
%


\begin{thebibliography}{99}

\bibitem{adler.ramanan} A. Adler and S. Ramanan, \emph{Moduli of
abelian varieties}, Lecture Notes in Math., {\bf 1644}, Springer,
Berlin, 1996.

\bibitem{bor.gunn.1}
L.~A. Borisov and P.~E. Gunnells, \emph{Toric varieties and modular forms},
Invent. Math., to appear, preprint {\tt math.NT/9908138}.

\bibitem{bor.gunn.2} L.~A. Borisov and P.~E. Gunnells, \emph{Toric
modular forms and nonvanishing of $L$-functions}, J. Reine. Angew. Math.,
to appear, preprint {\tt math.NT/9910141}.

\bibitem{bor.gunn.3} L.~A. Borisov and P.~E. Gunnells, \emph{Toric
modular forms of higher weight}, in preparation.

\bibitem{Chandra} K. Chandrasekharan, \emph{Elliptic functions}, Fundamental
Principles of Mathematical Sciences, {\bf 281}, Springer-Verlag, Berlin-New York, 1985.

\bibitem{dolgachev} I.~V. Dolgachev, \emph{Invariant stable bundles
over modular curves $X(p)$}, in {\it Recent progress in algebra
(Taejon/Seoul, 1997)}, 65--99, Contemp. Math., {\bf 224},
Amer. Math. Soc., Providence, RI, 1999.

\bibitem{edge} W.~L. Edge, \emph{Klein's encounter with the simple group of order
660},  Proc. London Math. Soc. {\bf 24} (1972), 647--668.

\bibitem{macaulay} D. Grayson and M. Stillman, \emph{Macaulay 2: A
computer program designed to support computations in algebraic
geometry and computer algebra.}  Source and object code available from
{\tt http://www.math.uiuc.edu/Macaulay2/}.

\bibitem{gross.popescu.1} M. Gross and S. Popescu, \emph{Equations of
$(1,d)$-polarized abelian surfaces}, Math. Ann. {\bf 310} (1998),
no.~2, 333--377.

\bibitem{gross.popescu.2} M. Gross and S. Popescu, \emph{The moduli
space of $(1,11)$-polarized abelian surfaces is unirational},
Compositio Math. (2000), to appear, preprint {\tt math.AG/9902017}.


\bibitem{Klein1} F. Klein, \emph{\"Uber transformationen siebenter
Ordnung der elliptischen Funktionen} (1878/79), Abhandlung {\bf
LXXXIV}, in {\it Gesammelte Werke}, Bd. {\bf III}, Springer, Berlin
1924.

\bibitem{Klein2} F. Klein, \emph{\"Uber die Transformation elfter
Ordnung der elliptischen Funktionen}, Math.  Ann. {\bf 15} (1879),
reprinted in {\it Ges. Math. Abh., Bd. {\bf III}, art.  {\bf LXXXVI},
pp. 140-168}.




\bibitem{Lang} S. Lang, \emph{Introduction to modular forms},
Springer-Verlag, 1976.

\bibitem{Mumford} D. Mumford, \emph{Tata lectures on theta I},  with the
assistance of C. Musili, M. Nori, E. Previato and M. Stillman. Progress in
Mathematics, {\bf 28}, Birkh\"{a}user Boston, Inc., Boston, Mass., 1983.

\bibitem{pol} A.~Polishchuk, \emph{Classical Yang-Baxter equation
and the $A_\infty$-constraint},\\ preprint {\tt math.AG/0008156}.

\bibitem{shimura} G. Shimura, \emph{Introduction to the arithmetic
theory of automorphic functions}, Reprint of the 1971 original,
Princeton Univ. Press, Princeton, NJ, 1994.

\bibitem{stark} H.~M. Stark, \emph{
Class fields and modular forms of weight one,} in {\it Modular functions of one variable, 
V (Proc. Second Internat. Conf., Univ. Bonn, Bonn, 1976)}, 277--287. 
Lecture Notes in Math., {\bf 601}, Springer, Berlin, 1977.

\bibitem{velu} J. V\'elu, \emph{Courbes elliptiques munies d'un
sous-groupe $\ZZ/n\ZZ\times {\mu}_{n}$.}  Bull. Soc. Math. France
M\'em. No. {\bf 57}, (1978), 5--152.

\end{thebibliography}
\end{document}